\documentstyle[12pt]{article}
\begin{document}
\title{ON COMMUTATION RELATIONS FOR QUONS
\footnote
{The paper is partially supported by KBN, Grant No 2P 302 087 06}}
\author{W{\l}adys{\l}aw Marcinek
\\ Institute of Theoretical Physics, University
of Wroc{\l}aw,\\ Pl. Maxa Borna 9, 50-204  Wroc{\l}aw,\\
Poland}
\date{}
\maketitle
\begin{abstract}
\noindent
The model of generalized quons is described in a purely
algebraic way. Commutation relations and corresponding
consistency conditions for our generalized quons system
are studied in terms of quantum Weyl algebras. Fock
space representation and corresponding scalar product
is also given.
\end{abstract}
\def\ot{\otimes}
\def\ra{\longrightarrow}
\def\ca{{\cal A}}
\def\cl{{\cal L}}
\def\cw{{\cal W}}
\def\cc{{\cal C}}
\def\cb{{\cal B}}
\def\br{[.,.]_{\Psi}}
\def\pr{(.|.)}
\def\sc{<.|.>}
\def\dr{der(\ca)}
\def\en{End(\ca)}
\def\sh{S_{\chi}}
\def\ps#1,#2,{\Psi_{{#1}{,}{#2}}}
\def\al{\alpha}
\def\ga{\gamma}
\def\ep{\epsilon}
\def\vp{\varphi}
\def\id{id_{\ca \ot \ca}}
\def\sa{S_{\ca,\ca}}
\def\1n{^{(1)}}
\def\2n{^{(2)}}
\def\3n{^{(3)}}
\def\p{\partial}
\def\kp{k-\Psi-}
\def\ba{\begin{array}}
\def\ea{\end{array}}
\def\be{\begin{equation}}
\def\ee{\end{equation}}
\newpage
\section{Introduction}
It is well known that all standard particles in physics like
electrons, photons, protons and all others can be divided into
two classes according to their statistics: fermions and bosons.
In the last years particle excitations equipped with new kind
of statistics have been discovered. For example particles with
fractional spin and statistics interpolating between fermions
and bosons appears in the two dimensional field theories. We
have for these particles the following commutation relations
for creation and anihilation operators
$$
a^*_i \ a_j - q \ a_j \ a^*_i = \delta_{ij} {\bf 1},
$$
where the deformation parameters $q$ is real and
$-1 \underline{<} q \underline{<} 1$ \cite{owg,gre,moh}.
Observe that the above relations for $q=1$ become the
canonical commutation relations for bosons and for $q=-1$
become the canonical anticommutation relations for fermions.
Note that for $q \neq \underline{+}1$ we have the mixed
relations between anihilation and creation operators only,
there are no relations for operators of the same kind, i.e.
for creation operators or for anihilation operators \cite{fi}.
Particles described by the above commutation relations called
quons have been studied by several authors, by D. Zagier \cite{zag},
by Bo$\dot{z}$ejko and Speicher \cite{bs1,bs2}, Finkelstein
\cite{fin}
and others. A simple physical model for particles with interpolating
statistics called anyons has been given by Wilczek, \cite{wil}.
According to his concept an anyon is a composite system which
consist a particle with charge $e$ and a singular magnetic field
concentrated completely on a vertical line. The charge is moving
around the singular this line in the region in which there is
no influence of the magnetic field. Note that every rotation
yields an arbitrary phase factor in the wave function of the
particle. This means that the statistics is fractional.
Anyons in general has been studied from different point of
view. The path integral method have been used by Wilczek, Zee and
Wu \cite{wil,wiz,swu}. The connection with the Chern-Simon
theory has been studied in Ref. \cite{ywu,asw,kkk}. Anyons
as $q$-deformation of bosons has been considered in Ref.
\cite{bdm}. In Ref. \cite{smg} the anyonic excitations are
described by the Zamolodchikov-Faddeev algebra. On the other hand
it is known that recently the deformed commutation relations
for creation and anihilation operators (CAO) has been studied from
different points of views by several authors. For example the
deformation of bosons and fermions corresponding to quantum
groups $SU_q(2)$ has been given by Pusz and Woronowicz \cite{P,PW}.
The deformation corresponding to superparticles has been considered
by Chaichian, Kulisch and Lukierski \cite{ckl}. Quantum deformations
have benn also studied by Fairle and Zachos \cite{FZ}, Vokos \cite{V},
Arik \cite{A} and many others. The commutation relations for Hecke
braiding has been studied by Kempf \cite{Ke}. It is interesting
that from algebraic point of view that all deformations of
commutation
relations for CAO can be described in terms of the so-called
Wick algebras \cite{js}. Note that Wick algebras are some
special examples of quantum Weyl algebras \cite{bae}, see
also Refers \cite{bor,ed}. Observe that in Wick algebra
there are no relations between creation (or anihilation)
operators themselves but such relations are possible in
certain cases \cite{js}. Obviously all such relations should
be consistent. Hence we need some additional assumption for
such consistency. Such conditions look like the Wess-Zumino
consistency conditions on quantum plane \cite{WZ}. Wick
algebras with additional relations and the corresponding
consistency conditions have been studied by the author
and by Ra{\l}owski in Ref. \cite{RM,m10,ral,mad}. In this
paper we are going to study a generalized quons as a generalization
of Wilczek model of anyons. For our generalized quons we obtain the
following relations
$$
\ba{ll}
a^*_i \ a_i - q_i \ a_i \ a^*_i = {\bf 1},\;&
a^*_i \ a_j - b_{ij} \ a_j \ a^*_i = 0,\mbox{for}\; i\neq j\;\\
a_i \ a_j - b_{ij} \ a_j \ a_i = 0,\;&
a^*_i \ a^*_j - b_{ij} \ a_j^* \ a^*_i = 0,
\ea
$$
where $-1 \underline{<} q_i \underline{<} 1$ are real parameters,
and $\overline{b}_{ij} = b_{ji}$ are some coefficients.
The paper is organized as follows. In Section 2 as the mathematical
preliminaries we recall the concept of quantum Weyl algebras.
In Sec. 3 we describe the model for generalized quons.
By a system of $N$ generalized quons we understood rotational
excitations corresponding to ordinary particle moving on a plane
with $N$ singularities. We assume that for every singular point
there is a corresponding operator of (discrete) rotation.
We also assume that such operators generate an abelian group.
We give the commutation relations for a system of our generalized
quons in terms of a quantum Weyl algebra.
Next we give the solution for the consistency conditions
corresponding
to a system of $N$-quons. The category of quantum states for our
system is described in Section 4. The Fock space representation is
considered in Section 5. The scalar product is also given.
\section{Quantum Weyl algebras}
Let us give a short introduction to the concept of quantum Weyl
algebras \cite{ed,bor} for the description of quantum states of
our generalized quons. Let $E$ be a finite dimensional Hilbert space
equipped with a basis $x^i, i = 1,...,N = dim H$. The conjugated dual
space is denoted by $E^*$. Note that our notation is not covariant.
We assume that the pairing $\pr : E^* \ot E \ra \bf C$ and the
corresponding scalar product
is given by
\be
(u^* | v) \equiv <u | v> := \Sigma_i^N \overline{u}^i v^i,
\ee
where $u^* = \Sigma_i^N \overline{u}^i x^{*i}$ and
$v = \Sigma_i^N v^i x^i$.
We have here the following\\
\begin{em}
{\bf Definition A :}
A linear and invertible operator $C : E^* \ot E \ra E \ot E^*$
\be
\ba{c}
C(x^{*i} \ot x^j) = \Sigma_{k,*l} \ C^{*ij}_{k*l} \ x^k \ot x^{*l} ,
\label{ce}
\ea
\ee
for which we have
\be
\ba{ccc}
(C^{*ij}_{k*l})^* = \overline{C}^{*ji}_{l*k},&
\mbox{i.e.}& C^* = \overline{C}^T,
\ea
\ee
where $(C^T)^{*ij}_{k*l} = C^{*ji}_{l*k}$,
is said to be {\it a skew} or {\it a cross} operator.\\
\end{em}
It is a fundamental object for all our next constructions.
Let $C: E \ot E^* \ra E^* \ot E$ be an arbitrary skew
operator. We define {\it a partial dual} to $C$ as a
linear operator $\tilde{C} : E \ot E \ra E \ot E$ such
that we have the relation
\be
\ba{c}
(\pr \ot id_E) \circ (id_{E^*} \ot \tilde{C})
= (id_E \ot \pr) \circ (C \ot id_E).
\label{pd}
\ea
\ee
on $E^* \ot E \ot E$.
We have $(x^{*i}|x^j) = \delta^{*ij}$, hence we obtain
\be
\ba{c}
(\tilde{C})^{ij}_{kl} = C^{*ki}_{l*j}.
\ea
\ee
Let us assume that we have a second linear and invertible operator
$B : E \ot E \ra E \ot E$.
\be
\ba{ll}
B(x^i \ot x^j) = B^{ij}_{kl} \ x^k \ot x^l ,&
\;\;\; (B^*)^{ij}_{kl} = \overline{B}^{*j*i}_{*l*k} .
\ea
\ee
\begin{em}
{\bf Definition C} An algebra $\cw \equiv \cw(B, C)$
\be
\cw \equiv \cw(B, C) := T(E \oplus E^*)/I_{B, C},
\ee
where $I_{B, C} = I_C + I_B + I_{B^*}$ is an ideal in
$T(E \oplus E^*)$ is said to be a quantum Weyl algebra
if it is generated by the following relations
\be
\ba{l}
x^{*i} \ x^j = \delta^{ij} \ {\bf 1} +
C^{*ij}_{k*l} \ x^k \ x^{*l},\\
x^i \ x^j - B^{ij}_{kl} \ x^k \ x^l = 0,\\
x^{*i} \ x^{*j} - \overline{B}^{*i*j}_{*k*l} \ x^{*k} \ x^{*l} = 0,
\label{ad}
\ea
\ee
and equipped with consistency conditions
\be
\ba{c}
B\1n B\2n B\1n = B\2n B\1n B\2n,\\
B^{(1)}C\2nC\1n = C\2nC\1n B^{(2)},\\
(id_{E \ot E} + \tilde{C})(id_{E \ot E} - B) = 0.
\label{cd}
\ea
\ee
\end{em}
Let us introduce a linear operator $R_n$ by the
formula
\be
\ba{l}
R_n := id + \tilde{C}^{(1)} + \tilde{C}^{(1)} \tilde{C}^{(2)}
+...+ \tilde{C}^{(1)}...\tilde{C}^{(n-1)} .
\label{pro}
\ea
\ee
Note that properties of the operator $R_2 := id + \tilde{C}$
is essential for the algebra $\cw(B, C)$, see Refers
\cite{bs2,js}. If the kernel of $R_2$ is trivial then we have
relations between elements of Wick algebra and their conjugates
only. There are no relations between elements of the some kind
themselves. Such additional relations are possible if and only
if the kernel of $R_2$ is nontrivial. Observe that if the operator
$\tilde{C}$ has an eigenvalue $\lambda = -1$, then $R_2$ has
nontrivial kernel and vice versa.
\section{Generalized quon particles}
Let us consider a particle moving on a plane with $N$
singular points. The proper physical nature of such
singularities is not important for us. Only discrete
rotations of this particle around singularities are
essential for our model. The moving between singularities
is completely irrelevant. Hence it is not difficult to show that
the rotational excitations on singularities behave like
a system of $N$ independent particles. These rotational
excitations are called {\it generalized quons}. In this
Section we are going to describe the Weyl algebra for our
system of generalized quons.
We assume that the $i$-th singular point is characterized by
a real numbers $q_i$, $-1 \underline{<} q_i \underline{<} 1$,
$i=1,...,N$. We also assume that $\sigma_i$ is an operator for
discrete rotation around the $i$-th singular point. Next we
assume for simplicity that $\sigma_i$, $(i=1,...,N)$ generate
an abelian group $\Gamma$. It is natural to assume that the
statistics for our quons is determined by a commutation factor
factor in the sense of Scheunert \cite{sch} on the group $\Gamma$.
Recall that a commutation factor on an abelian group $\Gamma$
is a mapping $b : \Gamma \times \Gamma \ra {\bf C} - \{0\}$
such that
\be
\ba{l}
b(\al + \beta,\ga) = b(\al,\ga)b(\beta,\ga),
\;\;\;b(\al,\beta)b(\beta,\al) = 1
\ea
\ee
for every $\al, \beta, \ga \in \Gamma$. We assume here in
addition that $b(\al,\al) = 1$ for every $\al \in \Gamma$.

Now let us consider the quantum  Weyl algebra $\cw(B, C)$
corresponding for the system of $N$ generalized quons.
In order to do that we assume that the element $x^i$
describe the quantum state corresponding to the one left
rotation around the $i$-th singular points, i.e $x^i$
correspond to the generator $\sigma_i$. In other words
the space $E$ is graded by $\Gamma$, i.e. $x^i$ is a
homogeneous element with respect to the gradation of
grade $\sigma_i$. Similarly $x^{*i}$ correspond to the
one right rotation around the $i$-th singularity. For the
gradation we have here the relation
\be
\mbox{grade}(x^i) + \mbox{grade}(x^{*i}) = 0.
\ee
In this way the Hilbert space $E$ can be understood as the
space of quantum states for left rotations and $E^*$
- for the right rotations. We need two operators
$C$ and $B$ for the construction of quantum Weyl algebra.
We assume that the operator $C$ is given by
\be
\ba{c}
C(x^{*i} \ot x^j) := c_{ji} \ x^j \ot x^{*i}
\;\;\;\mbox{(no sum)}
\label{cop}
\ea
\ee
for $i, j = 1,...,N$, where $q_{ij}$ are complex
parameters such that $c_{ii} \equiv q_i$, where
$q_i$ characterize the $i$-th singular point on a plane.
For the additional operator $B$ we assume that
\be
\ba{c}
B(x^i \ot x^j) := b_{ij} \ x^j \ot x^i ,
\;\;\;\mbox{(no sum)}
\label{bop}
\ea
\ee
where $b_{ij} := b(\sigma_i,\sigma_j)$ are coefficient
for a commutation factor
$b : \Gamma \times \Gamma \ra {\bf C} - \{0\}$
on the grading group $\Gamma$. This means that
$b_{ij}b_{ji} = 1$ and $b_{ii} = 1$ for all $i=1,...,N$.
In other words $B$ is an even colour unital braiding,
see Ref. \cite{WM4} for more details. From consistency
conditions we obtain
\be
c_{ij} =
\left\{
\ba{c}
b_{ij} \;\; \mbox{for} \;\; i \neq j\\
q_i \;\; \mbox{for} \;\; i = j
\label{csy}
\ea
\right.
\ee
Note that such solution for consistency condition has
been found by Borowiec and Kharchenko, see Ref. \cite{bk}
for differential calculi.
In this way the quantum Weyl algebra $\cw(B, C)$ for quons is an
algebra generated by $x^i$ and $x^{*j}$ (i,j = 1,...,N)
subject to relations
\be
\ba{ll}
x^{*i} x^i = 1 + q_i \ x^{i} x^{*i}\;\;\mbox{for}\;\;i = j,&
x^{*i} x^j = b_{ji} \ x^{j} x^{*i}\;\;\mbox{for}\;\;i \neq j,\\
x^{i} x^j = b_{ij} \ x^{j} x^{i}\;\;\mbox{for all}\;\; i,j,&
x^{*j} x^{*i} = b_{ji} \ x^{*i} x^{*j}\;\;
\mbox{for all}\;\;i,j.
\ea
\ee
This algebra is denoted by $\cw_{q,b}(N)$. Let us consider
a few simple examples:\\
{\bf Example 2.1} In the case of $N = 1$ we obtain
\be
x^* \ x = q \ x \ x^* + 1,
\ee
i.e. the usual $q$-commutation relations, see Ref. \cite{fi}
for example. The corresponding quantum Weyl algebra is denoted
by $\cw_q$.\\
{\bf Example 2.2} For $N=2$ we have an algebra generated by
$x, y$ and $x^*, y^*$ and relations
\be
x^* \ x = q \ x \ x^* + 1,\;\;
y^* \ y = p \ y \ y^* + 1,\;\;
x \ y = \underline{+} \ y \ x ,
\ee
where $q_1 = q$, $q_2 = p$, and $b_{ij} = \underline{+}1$.\\
{\bf Example 2.3} For arbitrary $N \underline{>} 2$ we can
define an algebra $\cw_{q,\Theta}(N)$ for which we have
\be
b_{ij} = \exp \frac{2\pi\Theta_{ij}}{n},
\ee
where $\Theta_{ij} = - \Theta_{ji}$, $\Theta_{ij}$ is an
integer (mod n) for every $i,j=1,...,N$, $i \neq j$,
and $q = \{q_i | -1 \underline{<} q_i \underline{<} 1, i=1,...,N\}$.
The case $q_i = \underline{+}$ for every $i=1,...,N$ corresponds
to particles which are in usual called anyons.
\section{The category of states}
In this Section we would like to describe all quantum
states for the system of generalized quons. As was
indicated in the previous Section, spaces $E$ and $E^*$
correspond for the left and right single rotations around
singularities. For multiply rotations we need tensor
product of these spaces. In this way all states can
be described by a family of spaces. Such family form
a braided monoidal category. Hence let us briefly recall
the concept of braided category.
A monoidal category $\cc \equiv \cc(\ot, k)$ is shortly
speaking a category $\cc$ equipped with a
monoidal associative operation (a bifunctor)
$\ot :{\cal C}\times{\cal C}\ra {\cal C}$
which has a two-sided identity object $k$,
see Ref. \cite{ML,jst,jet,Maj,bm} for details.
A family of natural isomorphisms
\be
\ba{c}
\Psi \equiv \{\Psi_{U \ot W} : U \ot W \ra W \ot U \}
\label{bs1}
\ea
\ee
such that we have the following relations
\be
\ba{l}
\Psi_{U \ot V ,W} = \Psi_{U,W} \1n \circ \Psi_{V,W} \2n,\\
\Psi_{U,V \ot W} = \Psi_{U,W} \2n \circ \Psi_{U,V}\1n,
\label{hex}
\ea
\ee
is said to be a braiding or a braid symmetry on $\cc$.
The corresponding category is said to be a braided monoidal
category \cite{Maj,ML}. It is denoted by $\cc \equiv \cc(\ot, k,
\Psi)$.
If in addition
\be
\ba{c}
\Psi_{U,V}^2 = id,
\label{kwad}
\ea
\ee
for every objects $U, V \in \cc$, then $\Psi$
is said to be a symmetry or tensor symmetry
and the corresponding category $\cc$ is
called a symmetric monoidal or tensor category, see
Ref. \cite{man,Lub,GRR}.
If the condition (\ref{kwad}) holds, then we
use the symbol $S$ instead of $\Psi$.
The following result is important for us.

If we have two linear and invertible operators
$C : E^* \ot E \ra E \ot E^*$, and $B : E \ot E \ra E \ot E$
such that the consistency conditions (\ref{cd}) are satisfied,
then there is a braided monoidal category $\cc \equiv(\ot, k, \Psi)$
generated by $E$, $E^*$, $C$ and $B$. The category $\cc$ contains:
the underground field $k = {\bf C}$ of complex numbers, the vector
space $E$ and the dual $E^*$, all tensor products of $E$ and $E^*$,
all direct sums, and some quotients. We define braidings
$\Psi_{E,E}, \Psi_{E^*,E}, \Psi_{E,E^*}$, and
$\Psi_{E^*,E^*}$ by formulae
\be
\ba{ll}
\Psi_{E,E} \equiv B ,&
\Psi_{E^*,E} \equiv C ,\\
\Psi_{E,E^*} \equiv C^{-1} ,&
\Psi_{E^*,E^*} \equiv B^* ,
\label{bre}
\ea
\ee
where $(B^*)^{ij}_{kl} = B^{*l*k}_{*j*i}$.
All above braidings can be extended
to all tensor products of spaces $E$ and $E^*$ by
formulae (\ref{hex}). In this way we obtain
\be
\ba{l}
\Psi_{E, E^{\ot l}} := B^{(l)} \circ \ldots
\circ B\1n ,\\
\Psi_{E^{\ot k}, E^{\ot l}} := \Psi_{E, E^{\ot l}}
^{(1)} \circ \ldots \circ \Psi_{E, E^{\ot l}}^{(k)},\\
\Psi_{E^{*\ot k}, E^{\ot l}}, \\
:= C^{(l)} \circ \ldots \circ C^{(1)} \circ
\ldots \circ C^{(k+l-1)} \circ \ldots \circ
C^{(k-1)} \circ C^{(k+l)} \circ \ldots \circ
C^{(k)}.
\label{bra}
\ea
\ee
Let us consider in particular the braided monoidal category
corresponding to our quons system. In this case we obtain
a category of $\Gamma$-graded spaces denoted by
$\cc(\ot, q, b, \Gamma)$. In this case we can understood
objects of the category as all possible spaces of quantum states
for our quons. For example the tensor product $E^{\ot k}$ can
be understood as a space for all quantum states corresponding
to the $k$ rotations to the left. If we substitute formulae
(\ref{cop}) and (\ref{bop}) into relations (\ref{bre}), then
we obtain the braid symmetry for quons. For example we have
\be
\ba{l}
\Psi_{E, E^{\ot l}}(x^i \ot (x^{j_1} \ot...\ot x^{j_l}))\\
:= \Pi_{k=1}^l \ b_{ij_k} \ ((x^{j_1} \ot...\ot x^{j_l}) \ot x^i).
\ea
\ee
Obviously the braid symmetric states should be equivalent.

An algebra $\ca \equiv \ca_B (E) := TE/I_{B}$,
is said to be $B$-symmetric algebra over $E$ \cite{Ba},
if it is generated by the relation
\be
\ba{c}
m \circ (id_{E \ot E} - B)(x^i \ot x^j) = 0,
\label{ide}
\ea
\ee
where $m$ is the multiplication in the algebra, or equivalently
\be
\ba{c}
x^i \ x^j - B^{ij}_{kl} \ x^k \ x^l = 0.
\label{mul}
\ea
\ee
This means that the algebra $\ca$ for our quon system
is an algebra generated by $\{x_i | i = 1,...,N \}$ and relation
\be
\ba{c}
x^i \ x^j = b_{ij} \ x^j \ x^i.
\label{rc}
\ea
\ee
This algebra can be regarded as the algebra of all braid
equivalent quantum states for left rotations.
On the other hand we have a $B$-symmetric algebra $\ca^*$ over
the space $E^*$, it is an algebra defined as the quotient
$\ca^* \equiv \ca_B (E^*) := TE^*/I^*$, such we have the
following relations
$$
x^{*i} \ot x^{*j} - \overline{B}^{*i*j}_{*k*l} \ x^{*k} \ot x^{*l} =
0.
$$
Obviously $\ca$ and $\ca^*$ are quadratic algebras over $E$, and
$E^*$, respectively. They are graded
\be
\ba{c}
\ca = \bigoplus_{n\in N} \ca^{n},\;\ca^0\equiv k,
\;\ca^1 \equiv E ,\\
\ca^* = \bigoplus_{n\in N} \ca^{*n},\;\ca^{*0}\equiv k,
\;\ca^{*1} \equiv E^* .
\ea
\ee
We denote by $P : TE \ra \ca$ and $P^* :
TE^* \ra \ca^*$ the quotient mappings. We have
\be
\ba{l}
P_k \equiv id_{E^{\ot k}} + B \ot id_{E^{\ot k-2}}
+ (B \ot id_{E^{\ot k-2}}) \circ (id_E \ot B \ot id_{E^{\ot k-3}}) +
...,\\
P^*_k \equiv id_{E^{*\ot k}} + B^* \ot id_{E^{*\ot k-2}}
+ (B^* \ot id_{E^{*\ot k-2}}) \circ (id_{E^*}
\ot B^* \ot id_{E^{* k-3}}) + ...,
\ea
\ee
for the $k$-th homogeneous components of $P$ and $P^*$ respectively.

Let us define a right evaluation mapping
$ev_k : E^* \ot E^{\ot k} \ra E^{\ot k-1}$ by formulae
\be
\ba{l}
ev_1 \equiv \pr,\\
ev_{k} := ev_1 \ot id_{E^{\ot k-1}} +
(id_{E^*} \ot ev_{k-1}) \circ (\tilde{C} \ot id_{E^{\ot k-1}}).
\label{eva}
\ea
\ee
\begin{em}
{\bf Lemma B} The evaluation mapping $ev$
is braid commutative, i.e. we have
\be
\ba{c}
ev_{k+l} \left[
id_{E^*} \ot (id_{E^{\ot k}\ot E^{\ot l}} - \ps E^{\ot k}, E^{\ot l},)
\right] = 0.
\label{wr}
\ea
\ee
on $E^* \ot E^{\ot k} \ot E^{\ot l}$.\\
\end{em}
{\bf Proof} We prove for $k=l=1$ only.
We have
$$
\ba{l}
ev_2 \left[
id_{E^*} \ot (id_{E \ot E} - B)
\right]\\
= \left[
(ev_1 \ot id_{E}) + (id_E \ot ev_1) \circ (\tilde{C} \ot id_E)
\right]
\left[
id_{E^*} \ot (id_{E \ot E} - B)
\right]\\
= (ev_1 \ot id_{E})
\left[
id_{E^*} \ot ((id_{E \ot E} + \tilde{C}) \circ (id - B))
\right] = 0.
\ea
$$
\hfill $\Box$\\
It follows from the braid commutativity of
the evaluation mapping $ev$ that it can be
projected to the algebra $\ca$. Hence the
projected evaluation is a mapping $\al :
E^* \ot \ca \ra \ca$ given by the relation
\be
\ba{l}
\al_1 := ev_1 \equiv \pr,\\
\al_k(id_{E^*} \ot P_k) := P_{k-1} \circ ev_k.
\ea
\ee
For the evaluation $\al$ we have the following
deformed Leibniz rule
\be
\ba{c}
\al_{k+l} \circ (id_{E^*} \ot m)
= m \circ \left[
(\al_k \ot id_{\ca^l}) + (id_{\ca^k} \ot \al_l)
\circ (\Psi_{E^*,\ca^k} \ot id_{\ca^l})
\right],
\label{le}
\ea
\ee
defined on $E^* \ot \ca^k \ot \ca^l$. We also have the
relation
\be
\ba{c}
\al_{l+1} \circ (id_{E^*} \ot \al_l) \circ
\left[
id_{E^* \ot E^* \ot \ca}
- (\Psi_{E^*,E^*} \ot id_{\ca})
\right] = 0
\label{leb}
\ea
\ee
on $E^* \ot E^* \ot \ca$.
We have for example
\be
\al_2 (x^{*i} \ot x^k x^l) = \delta^{ik} \ x^l
+ b_{ik} \ \delta^{il} \ x^k .
\ee
\section{Fock space representation}
We going to construct a representation $a : \cw(B, C) \ra End(\ca)$
of the Wick algebra $\cw$ on the algebra $\ca$. For generators $x^i$
and $x^{*j}$, $(i, j = 1,...,N)$ we define
\be
\ba{c}
a_{x^i} f \equiv m(x^i \ot f),
\;\;a_{x^{*i}} f \equiv \al(x^{*i} \ot f)
\label{cao}
\ea
\ee
for every $f \in \ca$.
One can extend this representation for the whole Wick algebra
$\cw$. Obviously our commutation relations for creation and
anihilation operators should be consistent with defining relations
for the algebra $\cw$. One can also consider the representation
of the Wick algebra $\cw$ on the algebra $\ca^*$ in a similar way.

Let us define a $\Psi$-bracket by the relation
\be
\ba{l}
[\hat{a}_i,\hat{a}_j]_{\Psi}
\equiv \br (\hat{a}_i \ot \hat{a}_j)
:= c \circ(id - \Psi)(\hat{a}_i, \ot \hat{a}_j),
\label{pbr}
\ea
\ee
where $\hat{a}_i$ stands for $a_{x^{*i}}$ or $a_{x^i}$, $c$
is the composition map and
\be
\ba{c}
\Psi(\hat{a}_{x^{*i}} \ot \hat{a}_{x^j}) :=
a_{\Psi(\hat{x}^{i} \ot \hat{x}^j)},
\ea
\ee
where $\hat{x}^i$ stands for $x^{*i}$ or $x^i$ and the
braiding $\Psi$ on the right hand of the above formula
is given by the relations (\ref{bre}).
We have here the following important\\
\begin{em}
{\bf Theorem A}
We have on the algebra $\ca$ the following commutation
relations for the representation of Wick algebra $\cw(B, C)$
\be
\ba{ccc}
[a_{x^{*i}}, a_{x^j}]_{\Psi} = \delta_{ij},&
[a_{x^i}, a_{x^j}]_{\Psi} = 0,&
[a_{x^{*i}}, a_{x^{*j}}]_{\Psi} = 0.
\label{cr2}
\ea
\ee
\end{em}
{\bf Proof} Using the relations (\ref{cao}) we obtain
\be
\ba{l}
(a_{x^{*i}} \circ a_{x^j})f \equiv
\al_{l+1} \circ (id_{E^*} \ot m) (x^{*i} \ot x^j \ot f),\\
(a_{x^i} \circ a_{x^j})f \equiv
m \circ (id_{E} \ot m)(x^i \ot x^j \ot f),\\
(a_{x^{*i}} \circ a_{x^{*j}}) f \equiv
\al_{l+1} \circ (id_{E^*} \ot \al_l) (x^{*i} \ot x^{*j} \ot f),
\ea
\ee
for $f \in \ca^l$. We also have
\be
\ba{l}
c \circ \Psi(a_{x^{*i}} \ot a_{x^j})f \equiv
m \circ (id_{E} \ot \al_{l+1}) \circ (C \ot id_{\ca})
(x^{*i} \ot x_j \ot f),\\
c \circ \Psi(a_{x^i} \ot a_{x^j})f \equiv
m \circ (id_{E} \ot m) \circ (B \ot id_{\ca} )(x^i \ot x^j \ot f),\\
c \circ \Psi(a_{x^{*i}} \ot a_{x^{*j}}) f \equiv
\al_{l+1} \circ (id_{E^*} \ot \al_{l})(B^* \ot id_{\ca})
(x^{*i} \ot x^{*j} \ot f),
\ea
\ee
It follows immediately from the definition
(\ref{pbr}) of the $\Psi$-bracket that we have
\be
[\hat{a}_i,\hat{a}_j]_{\Psi} f
= \left[
\hat{a}_i \circ \hat{a}_j - c \circ \Psi(\hat{a}_i, \ot \hat{a}_j)
\right]f.
\ee
Now we can calculate the first relation (\ref{cr2}) as follows
\be
\ba{l}
[a_{x^*i}, a_{x^j}]_{\Psi} f\\
= \left[
\al_{l+1} \circ (id_{E^*} \ot m)(x^{*i} \ot x^j \ot f)
- m \circ (id_{E} \ot \al_l) \circ (C \ot id_{\ca})
\right]\\
(x^{*i} \ot x^j \ot f)
= m \circ (\al_{1} \ot id_{\ca})(x^{*i} \ot x^j \ot f)
= \delta^{ij} f ,
\ea
\ee
where $f \in \ca$ and the relation (\ref{le}) has been used.
For the the second relation (\ref{cr2}) we obtain
\be
\ba{l}
[a_{x^i}, a_{x^j}]_{\Psi} f\\
= m \circ (id_E \ot m) \left[
id_{E \ot E \ot \ca} - (B \ot id_{\ca})(x^i \ot x^j \ot f)
\right] = 0.
\ea
\ee
Finally we calculate the third commutation relation (\ref{cr2})
\be
\ba{l}
[a_{x^{*i}}, a_{x^{*j}}]_{\Psi} f\\
 = \al_{l+1} \circ (id_{E^*} \ot \al_l)
\left[
id_{E^* \ot E^* \ot \ca} - (B^* \ot id_{\ca})
\right]
(x^{*i} \ot x^{*j} \ot f) = 0,
\ea
\ee
where the relation (\ref{leb}) has been used.
\hfill $\Box$\\
We introduce here the following state vector
\be
|n_1...n_N> = (x^1)^{n_1} ... (x^N)^{n_N}
\ee
up to some normalization factor. For the ground state vector
$|0> \equiv {\bf 1}_{\ca}$ we have as usual $a_{x^*i}|0> = 0$.
The scalar product can be given by the following formula
\be
\ba{l}
<s|t>_C := \delta_{mn} \ <s|P'_n t>
\label{sca}
\ea
\ee
for $s \in \ca^m$ and $t \in \ca^n$, where $P'_n$ is an operator
defined by induction
\be
P_{n+1} := (id \ot P_n) \circ R_n,
\ee
where $P_1 \equiv id$ and the operator $R_n$ is given by the
formula (\ref{pro}). Note that there is a very interesting
result of Bo$\dot{z}$ejko and Speicher \cite{bs2}. If the operator
$\tilde{C}$ is a bounded operator acting on some Hilbert space
such that
\be
\ba{ll}
(i)&\tilde{C} = \tilde{C}^* ,\\
(ii)&||\tilde{C}|| \underline{<} 1,\\
(iii)&\tilde{C}\;\;\mbox{is Yang-Baxter},
\ea
\ee
then the scalar product given by the relation (\ref{sca}) is positive
definite. Let us consider the representation of the Wick algebra
$\cw_{q,b}(N)$. corresponding for quons. In this case we have the
following action for creation and anihilation operators
\be
\ba{l}
a_{x^i} |n_1...n_N> = \Pi_{k=1}^{i-1} \ b_{ik} \
|n_1,...,n_i + 1,...,n_N>,\\
a_{x^{*i}} |n_1...n_N> = \Pi_{k=1}^{i-1} \ b_{ik} \
[n_i]_{q_i} |n_1...n_i -1... n_N>,
\ea
\ee
$[n_i]_{q_i} = \frac{1-q_i^{n_i}}{1-q_i}$, $(i=1,...,N)$.
For commutation relations we obtain
\be
a^*_i \ a_j - c_{ij} \ a_j \ a^*_i = \delta_{ij} {\bf 1},\;
a_i \ a_j - b_{ij} \ a_j \ a_i = 0,\;
a^*_i \ a^*_j - b_{ij} \ a_j^* \ a^*_i = 0,
\ee
where $a^*_i := a_{x^{*i}}$ and $a_i := a_{x^i}$.
It follows from the mentioned above theorem of Bo$\dot{z}$ejko
and Speicher that the corresponding scalar product is
positive definite if all $q_i$ are real,
$-1 \underline{<} q_i \underline{<} 1$, and
$\overline{b}_{ij} = b_{ij}$.
Note that in our quonic interpretation the vector $|n_1...n_N>$
describe the state corresponding to $n_1$ left rotations corresponding
to the first singularity,...,and to $n_N$ left rotations corresponding
to the last singularity. It is interesting that creation and
anihilation operators create and anihilate rotations (or more
precisely - rotational excitations) not particles! Observe that
in the above representation we have only left rotations. For
right rotations we need the representation on the algebra $\ca^*$.\\
\vspace{1cm}

\noindent {\bf Acknowledgments}\\
The author would like to thank to J.\ Lukierski and A.\ Borowiec
for discussion and critical remarks, to C.\ Juszczak, and R.\
Ra{\l}owski for any other help.

\end{document}